\title{Appendix to ``Real second-order freeness and the asymptotic real second-order freeness of several real matrix models'': Examples and Diagrams}
\author{C.\ Emily I.\ Redelmeier\thanks{Research supported by a two-year Sophie Germain post-doctoral scholarship provided by the Fondation math\'{e}matique Jacques Hadarmard, held at the D\'{e}partement de Math\'{e}matiques, UMR 8628 Universit\'{e} Paris-Sud 11-CNRS, B\^{a}timent 425, Facult\'{e} des Sciences d'Orsay, Universit\'{e} Paris-Sud 11, F-91405 Orsay Cedex.  E-mail address: emily.redelmeier@math.u-psud.fr}}

\documentclass[11pt]{article}

\usepackage{graphicx}
\usepackage{graphics}
\usepackage{amsfonts}
\usepackage{amssymb}
\usepackage{amsthm}
\usepackage{amsmath}
\usepackage{color}
\usepackage{cite}

\newtheorem{theorem}{Theorem}[section]

\theoremstyle{remark}

\newtheorem{example}[theorem]{Example}

\theoremstyle{definition}

\begin{document}

\maketitle

\begin{abstract}
We present examples and diagrams illustrating the proofs appearing in \cite{2011arXiv1101.0422R}, to which this paper is meant to be an appendix.  We show how matrix calculations may be represented by topological surface gluings, which may in turn be represented as permutations.
\end{abstract}

\section{Introduction}

Random matrix calculations can often be represented by topological constructions (see \cite{MR1492512, MR2036721, MR2052516}).  These topological constructions can be represented by permutations, as in \cite{MR0404045, MR746795, MR2052516}, which are often more convenient for the proofs.  This document is intended as an appendix to \cite{2011arXiv1101.0422R} (which we will refer to as the main paper), in which we attempt to provide some intuition for how the topological constructions are represented in the permutations, and describe the pictures which motivated the various proofs.  We illustrate with a number of examples.

In Section~\ref{section: face gluings}, we show how a matrix calculation involving Gaussian random matrices can be interpreted in terms of face gluings.  As pointed out in the exercises to \cite{MR2036721}, Chapter 3, calculations with real matrices include nonorientable surfaces.  In Section~\ref{section: cartography} we provide some intuition for how the permutations, and in particular premaps, may be used to represent these face gluings, even in the nonorientable case, and thus used in the matrix calculations.  In Section~\ref{section: zoo} we suggest pictures for the Ginibre and GOE matrices (Wishart matrices are discussed in the main paper).  We suggest intuition for further quantities calculated for the matrices: in Section~\ref{section: cumulants} we describe how cumulants are given by connected surfaces, in Section~\ref{section: PIE} how expressions with centred terms are given by diagrams where the regions corresponding to those terms are not disconnected and in Section~\ref{section: genus} we descibe asymptotic quantities in terms of topology and noncrossing diagrams.  In Section~\ref{section: freeness}, we describe how the conditions discussed in Section~\ref{section: PIE} and Section~\ref{section: genus} constrain the possible diagrams contributing to the cumulants of traces of products of centred matrices, those considered in first- and second-order freeness.

\section{Matrix Calculations as Face Gluings}
\label{section: face gluings}

As an example of how a topological surface may be used in a matrix calculation, let \(f_{ij}\) be independent \(N\left(0,1\right)\) random variables, and let \(X\) be an \(M\times N\) matrix with \(X_{ij}=\frac{1}{\sqrt{N}}f_{ij}\).  Let \(Y_{1},\ldots,Y_{n}\) be random matrices independent from \(X\), and let \(\gamma\in S_{n}\) and \(\varepsilon:\left[n\right]\rightarrow\left\{1,-1\right\}\).  In Lemma~3.4 of the main paper, we calculate expressions of the form
\[\mathbb{E}\left(\mathrm{tr}_{\gamma}\left(X^{\left(\varepsilon\left(1\right)\right)}Y_{1},\ldots,X^{\left(\varepsilon\left(n\right)\right)}Y_{n}\right)\right)\textrm{.}\]
We outline the topological constructions corresponding to the calculations in the proof.

\begin{figure}
\centering
\scalebox{0.5}{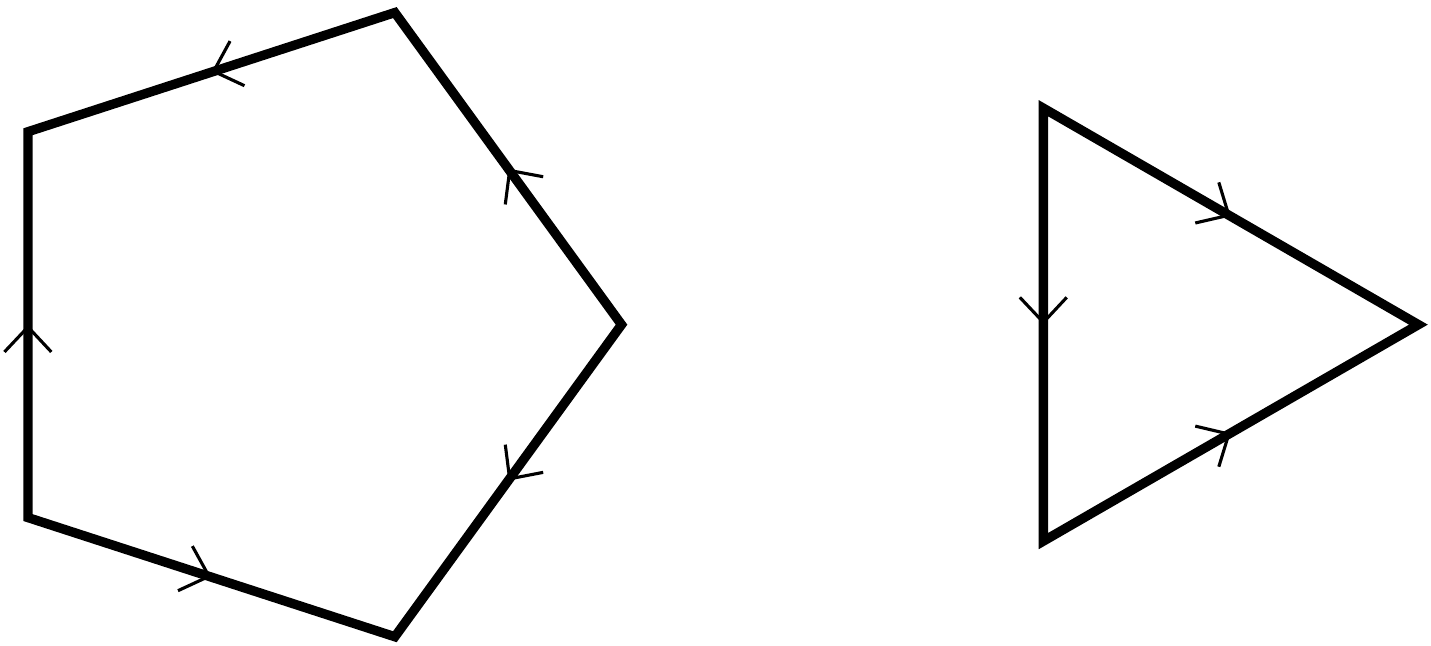}
\caption{The faces used to calculate the expected value of a product of traces of products of Gaussian random matrices given in the example.}
\label{figure: faces}
\end{figure}

Lemma~3.1 of the main paper suggests that the constraints on the indices appearing in a trace have a cyclic structure.  We construct a face for each trace in the expression, with the edges representing the terms from \(X\) and \(X^{T}\), and the vertices the \(Y_{k}\), cyclically in the order they appear in the cycles of \(\gamma\).  In Figure~\ref{figure: faces} we show the faces we would construct to calculate
\[\mathbb{E}\left(\mathrm{tr}\left(XY_{1}XY_{2}X^{T}Y_{3}Y_{4}X^{T}Y_{5}\right)\mathrm{tr}\left(X^{T}Y_{6}XY_{7}XY_{8}\right)\right)\textrm{.}\]
(We choose indices so that the \(k\)th occurrence of the matrix \(X\) appears with indices \(\iota^{+}_{k}\) and \(\iota^{-}_{k}\), which means the order they appear in depends on whether this occurrence of \(X\) appears with a transpose or not.  The arrows designate this direction.  The shared indices are shown between the matrices sharing them.)  We apply the Wick formula (see \cite{MR2036721}, p. 164) to the random variables.  The entries of \(X\) are independent random variables, so for a given pairing, the indices on the paired random matrix entries must be equal, which results in another set of constraints on the indices.  We represent these constraints by gluing the paired edges together so that the indices which must be equal are next to each other.  In this way we form a surface.

\begin{figure}
\centering
\scalebox{0.5}{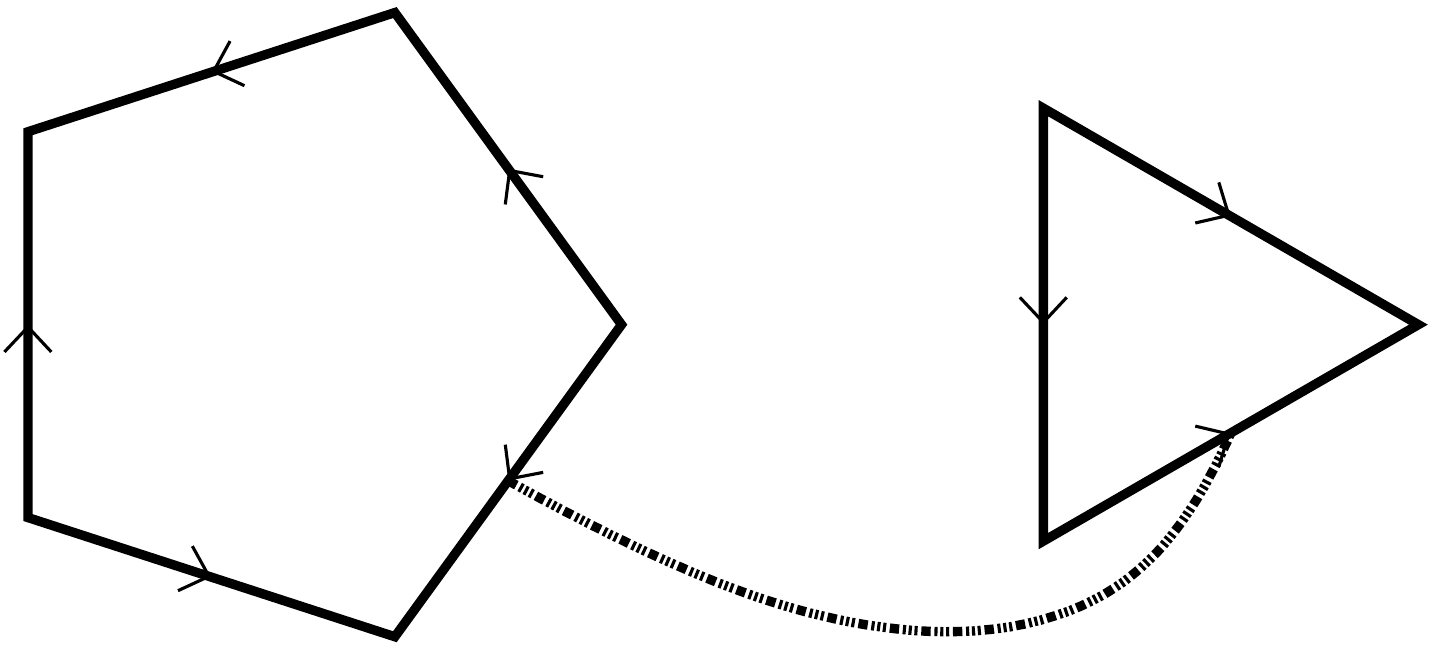}
\caption{An untwisted edge identification occurring when the term \(f_{\iota^{+}_{5}\iota^{-}_{5}}\) from the fifth random matrix term \(X^{T}\) is constrained to be equal to the term \(f_{\iota^{+}_{8}\iota^{-}_{8}}\) from the eighth matrix term \(X\).}
\label{figure: untwisted gluing}
\end{figure}

\begin{figure}
\centering
\scalebox{0.5}{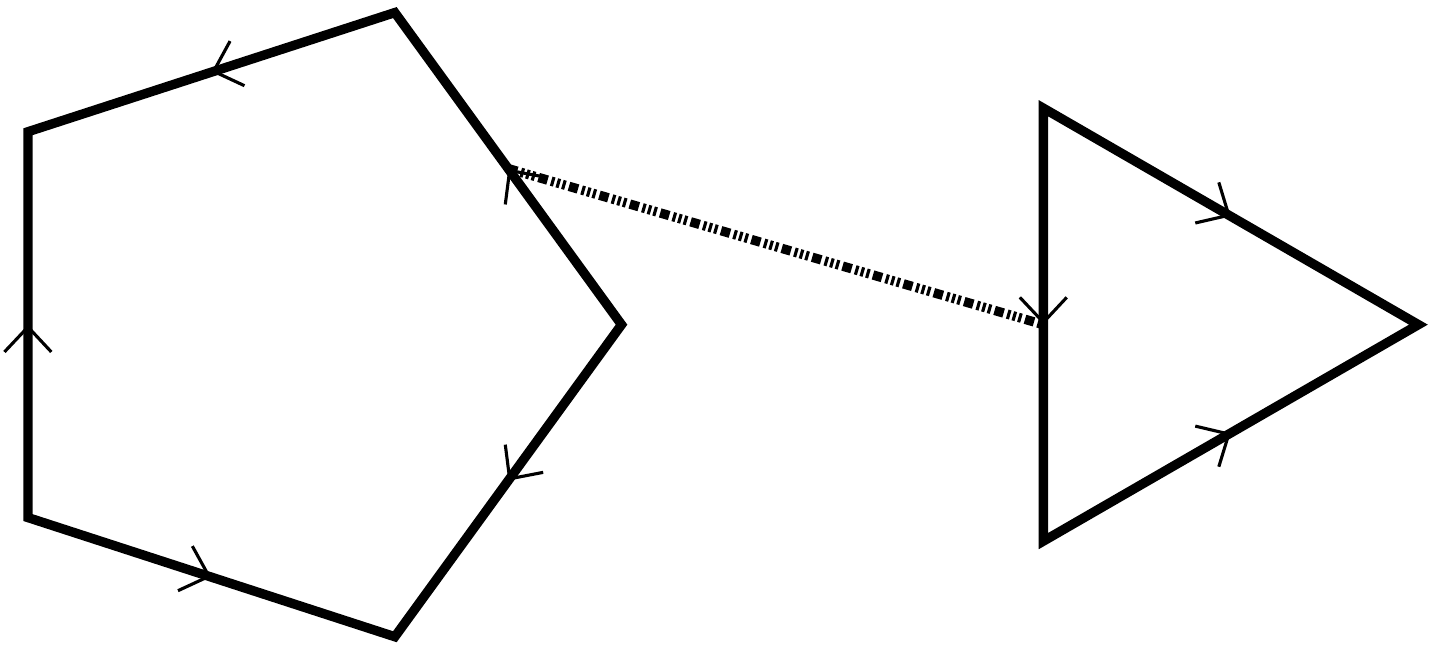}
\caption{A twisted edge identification occurring when the term \(f_{\iota^{+}_{1}\iota^{-}_{1}}\) from the first random matrix term \(X\) is constrained to be equal to the term \(f_{\iota^{+}_{8}\iota^{-}_{8}}\) from the seventh matrix term \(X\).}
\label{figure: twisted gluing}
\end{figure}

We note that if one term is from an \(X\) and another from \(X^{T}\) (as in Figure~\ref{figure: untwisted gluing}), the ordering of the indices means that the edge identification is untwisted.  However, if both are from an \(X\) or an \(X^{T}\) (as in Figure~\ref{figure: twisted gluing}), then the edge identification is twisted.  This is where differences from the complex case begin to appear.  If \(Z\) is a complex Gaussian random variable, then \(\mathbb{E}\left(Z^{2}\right)=\mathbb{E}\left(\overline{Z}^{2}\right)=0\), while \(\mathbb{E}\left(Z\overline{Z}\right)=1\).  Thus one term of a paired term must be from an \(X\) and the other from an \(X^{*}\), so the twisted edge identifications do not happen.  This excludes nonorientable surfaces, as well as the need to consider other relative orientations of the faces in diagrams drawn in the plane.  It is the latter that leads to the differences between the asymptotic fluctuations of the real and complex matrices.

\begin{figure}
\centering
\scalebox{0.5}{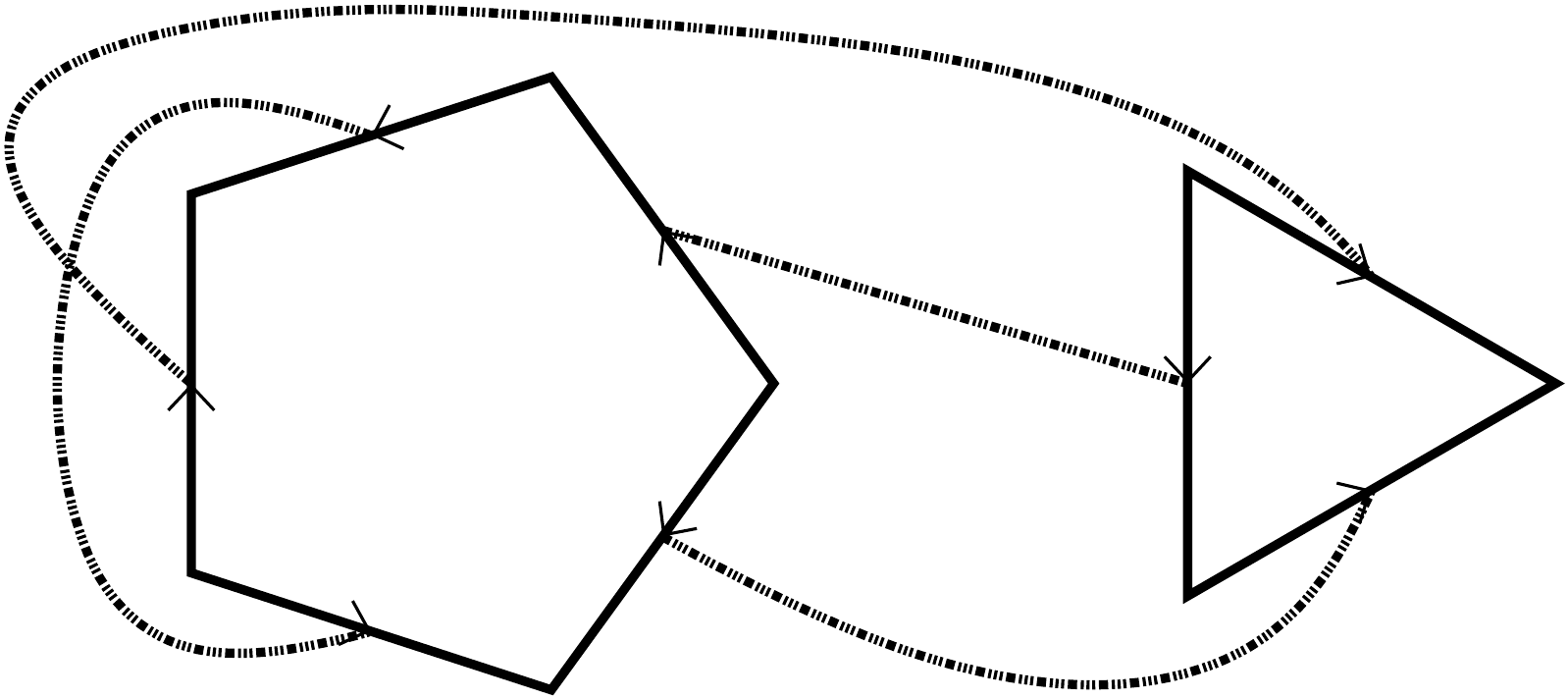}
\caption{An example of a gluing of the faces shown in Figure~\ref{figure: faces}.}
\label{figure: gluing}
\end{figure}

\begin{figure}
\centering
\scalebox{0.5}{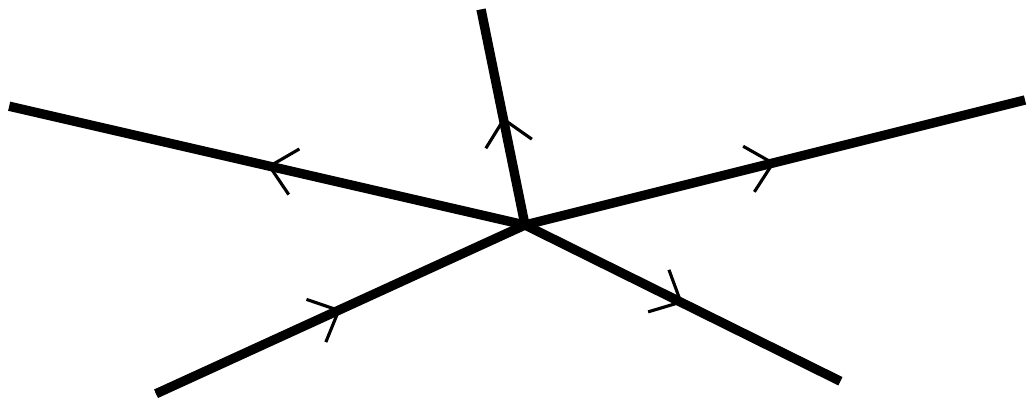}
\scalebox{0.5}{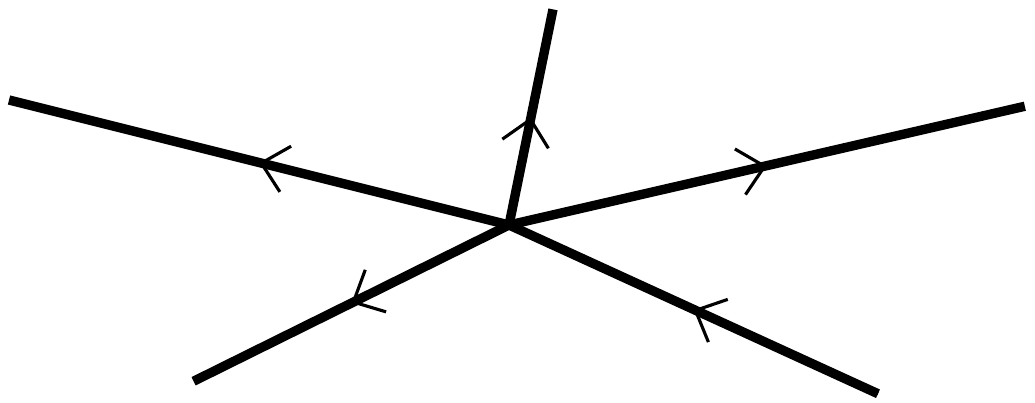}
\caption[Example of a vertex from a surface gluing]{A vertex which appears in the surface gluing shown in Figure~\ref{figure: gluing}, shown from both sides.}
\label{figure: vertex}
\end{figure}

We show a surface gluing in Figure~\ref{figure: gluing}.  The resulting surface will have a number of vertices.  The constraints on the indices around each are those that would appear in a trace, so each vertex corresponds to a trace of the \(Y_{k}\) matrices (see Figure~\ref{figure: vertex}).  We note that some of the corners of the faces may be upside down relative to the others.  Since this reverses the order of the indices as they appear, we replace that \(Y_{k}\) with \(Y_{k}^{T}\).  (In general, we can imagine that the transposes of the matrices are written on the backs of the faces.)  We note that if we consider the same vertex from the opposite side, the matrix product is the transpose of that on the original side, so the trace is unchanged.  In this case, the vertex contributes
\[\mathrm{Tr}\left(Y_{1}Y_{3}^{T}Y_{6}Y_{5}^{T}Y_{7}^{T}\right)=\mathrm{Tr}\left(Y_{7}Y_{5}Y_{6}^{T}Y_{3}Y_{1}^{T}\right)\textrm{.}\]

If we are considering an expression in normalized traces, a more natural quantity when we consider large \(N\) limits, then each trace contributes a factor of \(N\).  Thus, the highest order terms are those with the most vertices, that is, those with the largest Euler characteristics.

\section{Cartography on Unoriented Surfaces}
\label{section: cartography}

The face gluings (and hence the constraints on the indices) may be expressed in terms of permutations.  However, the encoding used in \cite{MR0404045, MR1603700, MR2036721} depends on the surface having an orientation.  We outline how we may similarly encode an unoriented face gluing.

If we have an oriented surface gluing, the ends of the edges have a cyclic ordering around the vertex they appear at.  We can define a permutation on the edge-ends \(\sigma\) (for {\em sommet}, French for vertex) with a cycle for each vertex in which the edge-ends at that vertex appear in their cyclic order, counterclockwise.  Similarly, we can define a permutation \(\alpha\) (for {\em ar\^{e}te}, edge) whose cycles contain the two ends of each edge.  Then the cycles of the permutation \(\varphi:=\sigma^{-1}\alpha^{-1}\) enumerate the edges appearing around each face of the surface gluing (with the convention that the edge is listed by its first edge-end encountered in a counterclockwise direction around the face; the other face sharing that edge will list the other edge-end).  A graph with this sort of surface embedding information, or this description of it in terms of permutations, is often called a map.

\begin{figure}
\centering
\scalebox{0.75}{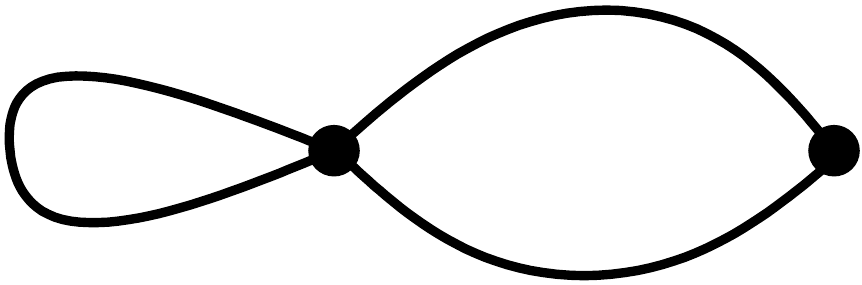}
\caption[A map]{A simple map, as described in Example~\ref{example: map}.}
\label{figure: map}
\end{figure}

\begin{example}
\label{example: map}
Consider the map in Figure~\ref{figure: map}.  The edge-ends are labelled by the integers in \(\left[6\right]\).  Then \(\sigma=\left(1,2,3,4\right)\left(5,6\right)\), and \(\alpha=\left(1,5\right)\left(2,3\right)\left(4,6\right)\).  We calculate that \(\sigma^{-1}\alpha^{-1}=\left(1,6,3\right)\left(2\right)\left(4,5\right)\).  The first cycle corresponds to the outside face, the second to the loop on the left, and the third to the two-edged face on the right.  We note that the face with two edges lists its edges as \(4\) and \(5\), while these edges are referred to as \(6\) and \(1\), respectively, when they are listed by the outside edge.
\end{example}

It is also possible to represent hypermaps in this way.  A hypermap is a map which may have hyperedges, that is, edges which may have any positive integer number of ends, rather than just two.  The ends of a hyperedge must also have a cyclic ordering representing their embedding in the surface.  Hyperedges may be visualized as either a separate set of vertices or faces, with the edge-ends becoming the edges.

\begin{figure}
\centering
\scalebox{0.75}{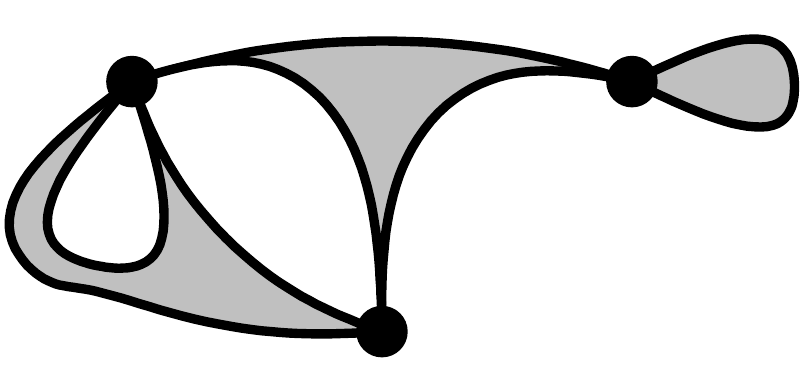}
\caption[A hypermap]{A hypermap map, as described in Example~\ref{example: hypermap}.}
\label{figure: hypermap}
\end{figure}

\begin{example}
\label{example: hypermap}
In Figure~\ref{figure: hypermap}, we show a hypermap.  The hyperedges are shown as outlined shapes connecting them to the black circles representing the vertices.  The edge-ends are labelled by integers in \(\left[7\right]\).  We have
\[\sigma=\left(1,2,3\right)\left(4,5\right)\left(6,7\right)\]
and
\[\alpha=\left(1,6,5\right)\left(2,7,3\right)\left(4\right)\textrm{.}\]
We can calculate that 
\[\sigma^{-1}\alpha^{-1}=\left(1,4,5,7\right)\left(2\right)\left(3,6\right)\textrm{.}\]
The first cycle corresponds to the outside face, the second to the region on the left with one vertex, and the third to the region with two vertices.
\end{example}

While both examples are planar in order to have clearer diagrams, it is possible to use these descriptions for maps on higher genus surfaces.  However, the description depends on the surface having a consistently defined orientation.  In order to describe the nonorientable surfaces which may appear in the real matrix problems, we extend this construction to the nonorientable case.

Any surface has an orientable covering space, which we may define by constructing two preimages of every point corresponding to the two possible orientations.  Points are close in the covering space if they are close in the original space with the same orientation.  For more details on this construction, see  \cite{MR1867354}, pages 234--235.  We can think of the covering space as the surface experienced by someone on the surface, rather than within it.

\begin{figure}
\centering
\scalebox{0.5}{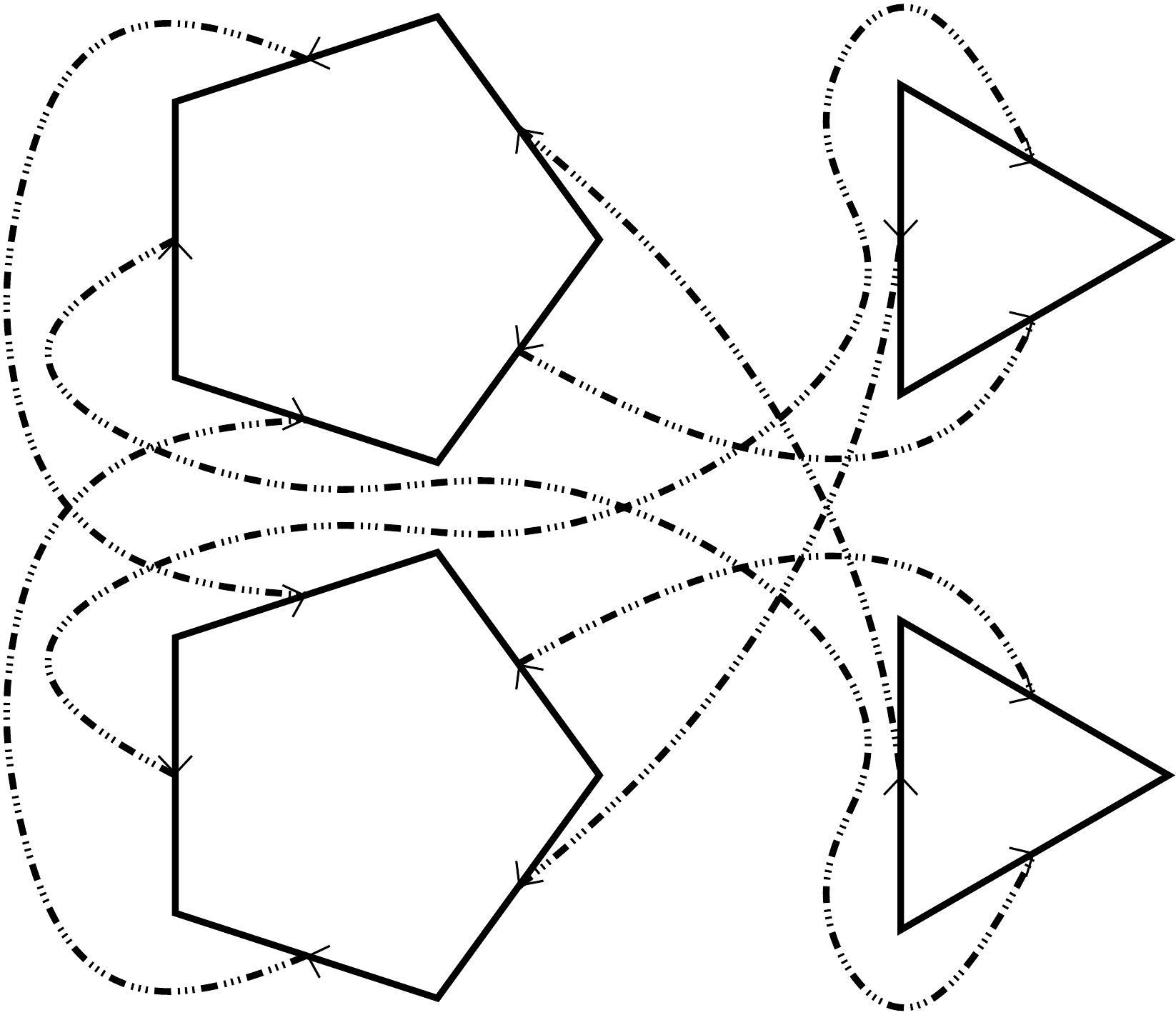}
\caption{The orientable covering space of the gluing shown in Figure~\ref{figure: gluing}.}
\label{figure: cover}
\end{figure}

For each face of our original surface, we construct two preimage faces, which we can think of as the front and the back.  We label the latter with negative signs.  When we glue edges in an untwisted way, we glue the corresponding edges on the fronts and the corresponding edges on the backs, while when we glue edges in a twisted way, we glue the corresponding edge of the front to the corresponding edge of the back, and vice versa.  If the original surface was nonorientable, it is possible to travel from the front to the corresponding point on the back within the covering space.  If the original surface was orientable, we obtain two disconnected copies of the original surface.  In Figure~\ref{figure: cover}, we show the cover constructed in this way from the gluing shown in Figure~\ref{figure: gluing}.

In the diagrams representing the matrix calculations in this work, we will generally begin with faces and edge or hyperedge information, and calculate vertex information, but this is simply the dual problem.  If we begin with face information described by a permutation \(\gamma\), we construct two preimage faces for each face: the fronts described by permutation \(\gamma_{+}\), and those on the back described by \(\gamma_{-}^{-1}\), since positive orientation is in the opposite direction on the back.  We note that \(\gamma_{+}\gamma_{-}^{-1}\), which represents all faces of the covering space, satisfies the definition for a premap.  In the example shown in Figure~\ref{figure: gluing}, we have
\[\gamma_{+}=\gamma=\left(1,2,3,4,5\right)\left(6,7,8\right)\textrm{,}\]
\[\gamma_{-}=\left(-1,-2,-3,-4,-5\right)\left(-6,-7,-8\right)\textrm{,}\]
and
\[\gamma_{+}\gamma_{-}^{-1}=\left(1,2,3,4,5\right)\left(6,7,8\right)\left(-5,-4,-3,-2,-1\right)\left(-8,-7,-6\right)\textrm{.}\]

The hyperedges may connect front and back, and therefore a cycle may contain both positive and negative integers.  However, for any cycle, we should be able to find the same cycle from the other side of the surface: a distinct cycle which contains the integers with sign reversed and in opposite order; thus, the permutation \(\pi\) representing the hyperedges should also be a premap.  See \cite{MR746795} for a similar construction.

We determine the vertex information in a similar manner.  Since the matrix \(Y_{k}\) appears counterclockwise of the \(k\)th occurrence of \(X\) or \(X^{T}\) on the front but clockwise of it on the back, the convention for which edge-end corresponds to which side of an edge means that the permutation \(\gamma_{-}\gamma_{+}^{-1}\pi^{-1}\) will list \(\gamma_{-}\left(-k\right)\) instead of \(-k\) for any negative integer.  We find that we must conjugate the permutation by \(\gamma_{-}^{-1}\), giving us the permutation \(\gamma_{+}^{-1}\pi^{-1}\gamma_{-}\).  The traces read off of the vertices are read clockwise, so we will be using the inverse of the permutation describing the vertices: \(\gamma_{-}^{-1}\pi\gamma_{+}\).

\begin{figure}
\centering
\scalebox{0.5}{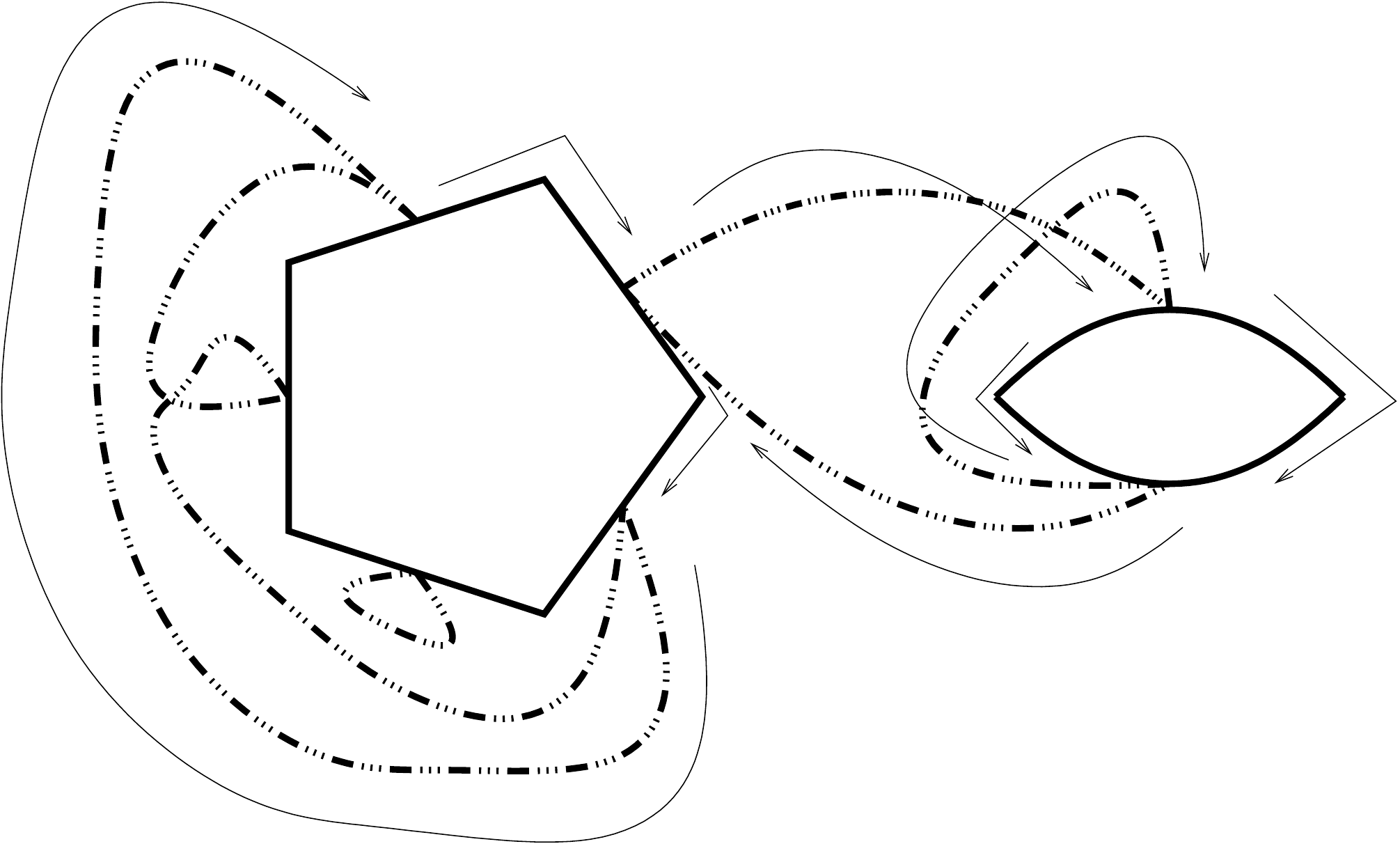}
\caption[A nonorientable hypermap]{An example of a nonorientable hypermap, discussed in Example~\ref{example: premap}.}
\label{figure: premap}
\end{figure}

\begin{example}
\label{example: premap}
Let \(\gamma=\left(1,2,3,4,5\right)\left(6,7\right)\) describe the faces of a nonorientable map, and let \(\pi=\left(1,7,-6\right)\left(6,-7,-1\right)\left(2,5,-3\right)\left(3,-5,-2\right)\left(4\right)\left(-4\right)\).  We show this construction in Figure~\ref{figure: premap}, and represent the change in sign in a cycle by a twisting of the hyperedge (this represents that the edge is glued in reverse direction).  We calculate that
\[\gamma_{+}^{-1}\pi^{-1}\gamma_{-}=\left(1,-6,7,5\right)\left(-5,-7,6,-1\right)\left(2,-3,-4\right)\left(4,3,-2\right)\textrm{.}\]
In the matrix calculations, we would use the inverse of this permutation.  We may trace the vertices in the diagram by following the boundaries of the faces and hyperedges, respecting twists.  The vertex \(\left(1,-6,7,5\right)\) is shown in the diagram.
\end{example}

\begin{figure}
\centering
\scalebox{0.75}{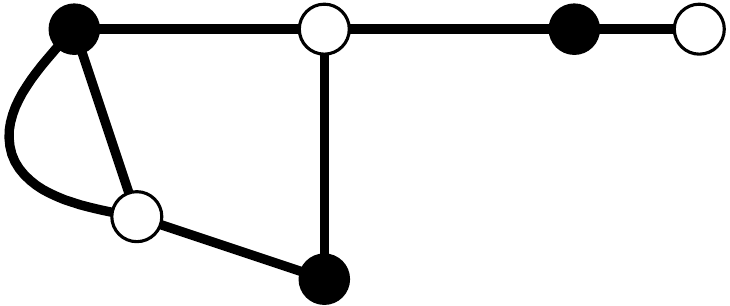}
\caption{The hypermap in Figure~\ref{figure: hypermap} where the hyperedges are seen as a set of vertices, shown in white, which alternate with the original vertices, and the edge-ends are seen as the edges.}
\label{figure: hyperedges as vertices}
\end{figure}

We may think of the hyperedges as either faces or vertices and the edge-ends as conventional edges.  In Figure~\ref{figure: hypermap}, hyperedges (in grey) may be seen as faces alternating with the original faces (in white), while the boundaries of the hyperedges (in black) may be seen as the edges (with the convention that each edge-end corresponds to the boundary clockwise from it relative to the hyperedge).  Figure~\ref{figure: hyperedges as vertices} shows the graph when we think of the hyperedges as vertices (white circles) which alternate with the original vertices (black circles), and again the edge-ends (black lines) are the edges.  In either case, we get the following formula for the Euler characteristic (Definition 3.6 in the main paper):
\[\chi\left(\gamma,\pi\right):=\#\left(\gamma_{+}\gamma_{-}^{-1}\right)/2+\#\left(\pi\right)/2+\#\left(\gamma_{-}^{-1}\pi^{-1}\gamma_{+}\right)/2-\left|I\right|\textrm{.}\]

\begin{example}
We return to the example shown in Figure~\ref{figure: gluing}.  We have
\[\pi=\left(1,-7\right)\left(-1,7\right)\left(2,-4\right)\left(-2,4\right)\left(3,-6\right)\left(-3,6\right)\left(5,8\right)\left(-5,-8\right)\textrm{.}\]
If the edge identification is untwisted, then in the covering space the front is glued to the front and back to the back, so the integers with the same sign appear in the same cycle.  If the edge identification is twisted, then the front is glued to the back and the back to the front, so integers with opposite signs appear in the same cycle.  We calculate that
\[\gamma_{-}^{-1}\pi\gamma_{+}=\left(1,-3,6,-5,-7\right)\left(7,5,-6,3,-1\right)\left(2,-8,-4\right)\left(4,8,-2\right)\]
and
\[\gamma_{-}^{-1}\pi\gamma_{+}/2=\left(1,-3,6,-5,-7\right)\left(2,-8,-4\right)\textrm{.}\]
The vertex represented by the first cycle is the one shown in Figure~\ref{figure: vertex}.

We calculate that \(\chi\left(\gamma,\pi\right)=2+4+2-8=0\).  Since the surface is connected and nonorientable, it must be a Klein bottle.
\end{example}

\section{The Matrix Models}
\label{section: zoo}

The expansion satisfied by all the matrix models in the main paper:
\begin{multline*}
\mathbb{E}\left(\mathrm{tr}_{\gamma}\left(X_{\lambda_{1}}^{\left(\varepsilon\left(1\right)\right)}Y_{1},\cdots,X_{\lambda_{n}}^{\left(\varepsilon\left(n\right)\right)}Y_{n}\right)\right)
\\=\sum_{\pi\in PM_{c}\left(\pm\left[n\right]\right)}N^{\chi\left(\gamma,\delta_{\varepsilon}\pi\delta_{\varepsilon}\right)-2\#\left(\gamma\right)}f_{c}\left(\pi\right)\mathbb{E}\left(\mathrm{tr}_{\gamma_{-}^{-1}\delta_{\varepsilon}\pi^{-1}\delta_{\varepsilon}\gamma_{+}/2}\left(Y_{1},\ldots,Y_{n}\right)\right)
\end{multline*}
may be interpreted as a set of rules for constructing surfaces out of the given faces.

In the case of Ginibre matrices, we may interpret the expression as a sum over all surface gluings which are consistent with directions built into the edges of the faces (counterclockwise if the matrix appears untransposed, and clockwise if it appears transposed).  We give an example:

\begin{figure}
\centering
\scalebox{0.75}{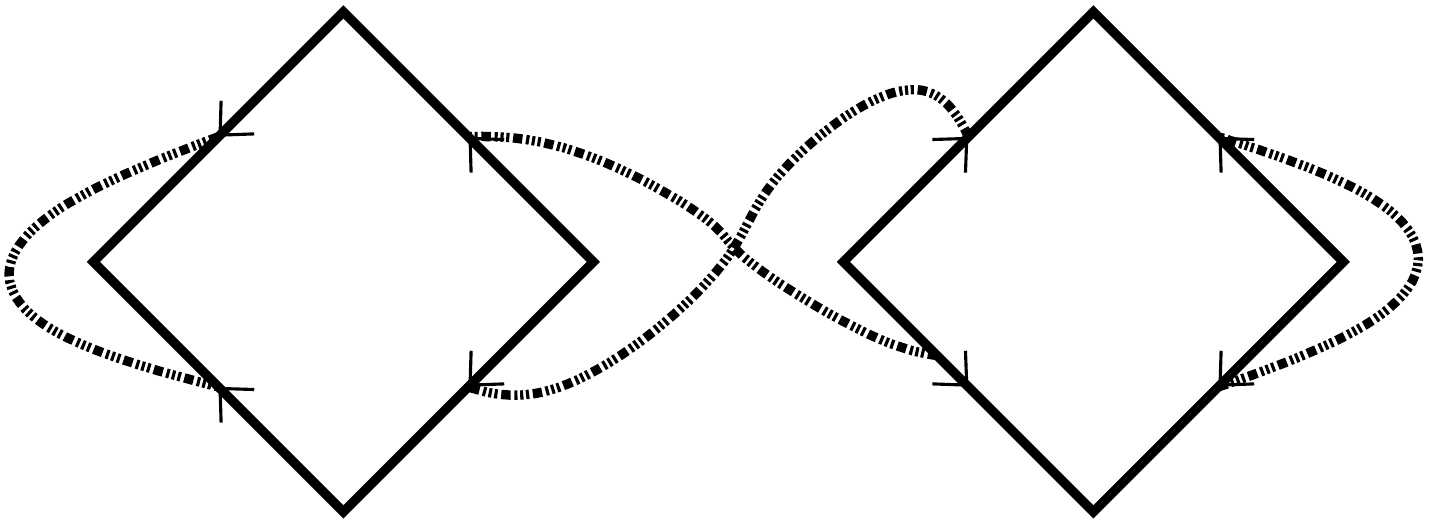}
\caption[Example of a face gluing for Ginibre matrices]{The faces constructed in order to compute the expression in Example~\ref{example: Ginibre matrices} with a possible edge identification marked.}
\label{figure: Ginibre matrices}
\end{figure}

\begin{example}
\label{example: Ginibre matrices}
Let \(Z\) be a real Ginibre matrix.  If we wish to calculate the expression
\[\mathbb{E}\left(\mathrm{tr}\left(ZY_{1}ZY_{2}Z^{T}Y_{3}Z^{T}Y_{4}\right)\mathrm{tr}\left(ZY_{5}Z^{T}Y_{6}ZY_{7}Z^{T}Y_{8}\right)\right)\textrm{,}\]
we construct the faces shown in Figure~\ref{figure: Ginibre matrices}.  We have \(\gamma=\left(1,2,3,4\right)\left(5,6,7,8\right)\) and \(\delta_{\varepsilon}=\left(1\right)\left(-1\right)\left(2\right)\left(-2\right)\left(3,-3\right)\left(4,-4\right)\left(5\right)\left(-5\right)\left(6,-6\right)\left(7\right)\left(-7\right)\left(8,-8\right)\).

The pairing shown is
\[\pi=\delta_{\varepsilon}\rho\delta\rho\delta_{\varepsilon}=\left(1,-7\right)\left(-1,7\right)\left(2,3\right)\left(-2,-3\right)\left(4,-6\right)\left(-4,6\right)\left(5,8\right)\left(-5,-8\right)\]
for \(\rho=\left(1,7\right)\left(2,3\right)\left(4,6\right)\left(5,8\right)\).  Paired integers of the same sign represent an untwisted edge identification and paired integers of opposite signs represent a twisted one.  We compute that
\[\gamma_{-}^{-1}\pi\gamma_{+}=\left(1,3,-5,-7\right)\left(7,5,-3,-1\right)\left(2\right)\left(-2\right)\left(4,-6\right)\left(6,-4\right)\left(8\right)\left(-8\right)\]
and that \(\chi\left(\gamma,\pi\right)=2+4+4-8=2\), so the example is a sphere (since it is connected), and contributes the term
\[\mathbb{E}\left(\mathrm{tr}\left(Y_{1}Y_{3}Y_{5}^{T}Y_{7}^{T}\right)\mathrm{tr}\left(Y_{2}\right)\mathrm{tr}\left(Y_{4}Y_{6}^{T}\right)\mathrm{tr}\left(Y_{8}\right)\right)N^{-2}\textrm{.}\]
\end{example}

The expression for GOE matrices may be interpreted as a sum over all edge gluings, in either direction.

\begin{figure}
\centering
\scalebox{0.75}{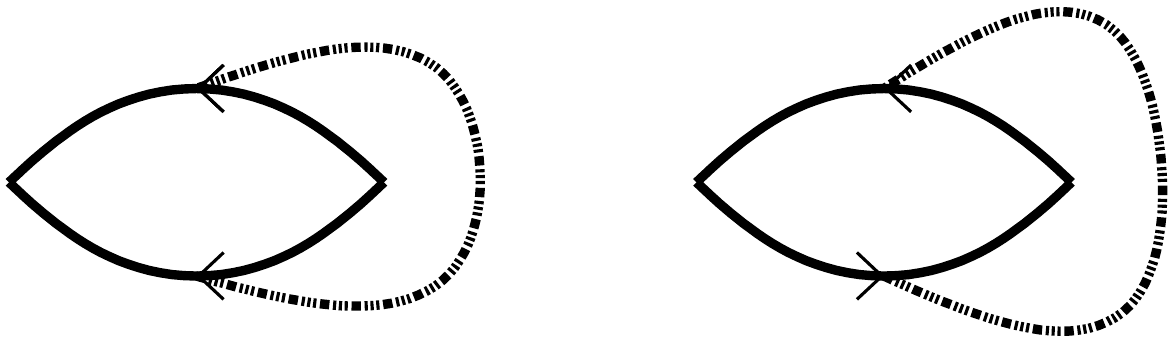}
\caption[Example of face gluings for GOE matrices]{The two surfaces corresponding to the two terms contributing to the expected value calculated in Example~\ref{example: GOE matrices}.}
\label{figure: GOE matrices}
\end{figure}

\begin{example}
\label{example: GOE matrices}
If we wish to calculate \(\mathbb{E}\left(\mathrm{tr}\left(T^{2}\right)\right)\), where \(T\) is a GOE matrix, we construct a face with two edges.  We may identify the two edges in an untwisted way, making the surface a sphere (Euler characteristic \(2\)), or in a twisted way, making the surface a projective plane (Euler characteristic \(1\)) (see Figure~\ref{figure: GOE matrices}).  These two surfaces correspond to the two elements of \(PM\left(\pm\left[2\right]\right)\cap{\cal P}_{2}\left(\pm\left[2\right]\right)\), that is, \(\left(1,2\right)\left(-1,-2\right)\) and \(\left(1,-2\right)\left(-1,2\right)\) respectively.  We calculate that
\[\mathbb{E}\left(\mathrm{tr}\left(T^{2}\right)\right)=1+N^{-1}\textrm{.}\]

If we wish to calculate \(\mathbb{E}\left(\mathrm{tr}\left(TY_{1}T_{2}Y_{2}\right)\right)\), where \(Y_{1}\) and \(Y_{2}\) are random matrices independent from \(T\), we get
\[\mathbb{E}\left(\mathrm{tr}\left(TY_{1}T_{2}Y_{2}\right)\right)=\mathbb{E}\left(Y_{1}\right)\left(Y_{2}\right)+\mathbb{E}\left(Y_{1}Y_{2}^{T}\right)N^{-1}\textrm{;}\]
the matrices \(Y_{1}\) and \(Y_{2}\) appear in different vertices in the sphere, but appear in the same vertex but opposite orientation in the projective plane.
\end{example}

An example calculation involving Wishart matrices and another involving several independent matrices are given in the main paper.

\section{Cumulants}
\label{section: cumulants}

\begin{figure}
\centering
\scalebox{0.5}{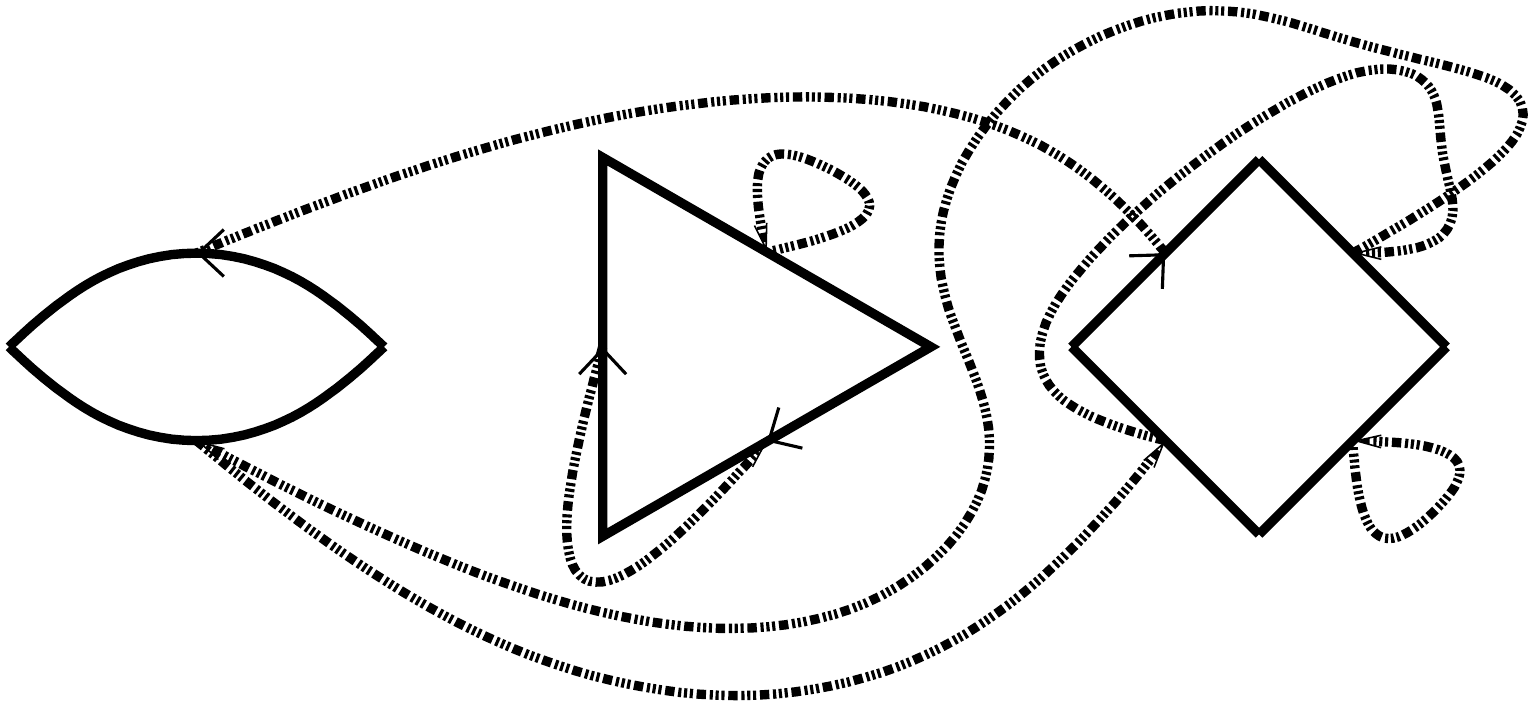}
\caption{A surface with two connected components.}
\label{figure: several matrices}
\end{figure}

Cumulants are sums over only connected surfaces: the moment-cumulant formula organizes the terms by which faces are part of a connected component.  The terms which appear in the term corresponding to partition (on elements corresonding to the traces, that is, the faces) \(\pi\) are those where the faces each block of \(\pi\) form a connected component.  The surface constructed in Figure~\ref{figure: several matrices} for the calculation of
\[\mathbb{E}\left(\mathrm{tr}\left(ZW_{2}^{\left(\lambda_{2}\right)}\right)\mathrm{tr}\left(W_{1}^{\left(\lambda_{3}\right)}Z^{T}Z^{T}\right)\mathrm{tr}\left(W_{2}^{\left(\lambda_{6}\right)}Z^{T}W_{2}^{\left(\lambda_{8}\right)}W_{1}^{\left(\lambda_{9}\right)}\right)\right)\]
would appear in the term
\[k_{2}\left(\mathrm{tr}\left(ZW_{2}^{\left(\lambda_{2}\right)}\right),\mathrm{tr}\left(W_{2}^{\left(\lambda_{6}\right)}Z^{T}W_{2}^{\left(\lambda_{8}\right)}W_{1}^{\left(\lambda_{9}\right)}\right)\right)k_{1}\left(\mathrm{tr}\left(W_{1}^{\left(\lambda_{3}\right)}Z^{T}Z^{T}\right)\right)\textrm{.}\]
These terms account for all the disconnected terms, so only the terms connecting all three faces, such as the one in Figure~\ref{figure: PIE}, contribute to the cumulant of the three traces:
\[k_{3}\left(\mathrm{tr}\left(ZW_{2}^{\left(\lambda_{2}\right)}\right),\mathrm{tr}\left(W_{1}^{\left(\lambda_{3}\right)}Z^{T}Z^{T}\right),\mathrm{tr}\left(W_{2}^{\left(\lambda_{6}\right)}Z^{T}W_{2}^{\left(\lambda_{8}\right)}W_{1}^{\left(\lambda_{9}\right)}\right)\right)\textrm{.}\]

\section{Centred Terms and the Principle of Inclusion and Exclusion}
\label{section: PIE}

When we expand a cumulant of traces of products of centred terms,
\begin{multline*}
k_{r}\left(\mathrm{tr}\left(\mathaccent"7017{A}_{1}\cdots \mathaccent"7017{A}_{p_{1}}\right),\ldots,\mathrm{tr}\left(\mathaccent"7017{A}_{p_{1}+\cdots+p_{r-1}+1}\cdots\mathaccent"7017{A}_{p}\right)\right)\\=\sum_{K\subseteq\left[p\right]}\left(-1\right)^{\left|K\right|}\prod_{k\in K}\mathbb{E}\left(\mathrm{tr}\left(A_{k}\right)\right)\\k_{r}\left(\mathrm{tr}\left(\prod_{k\in\left[1,p_{1}\right]\setminus K}A_{k}\right),\ldots,\mathrm{tr}\left(\prod_{k\in\left[p_{1}+\cdots+p_{r-1}+1,p\right]\setminus K}A_{k}\right)\right)\textrm{,}
\end{multline*}
we can think of the terms with indices in subset \(K\), which appear in the first product, as being required to have a diagram not connected to any other (that is, one that would contribute to \(\mathbb{E}\left(\mathrm{tr}\left(A_{k}\right)\right)\).  Other terms may also be disconnected.

\begin{figure}
\centering
\scalebox{0.5}{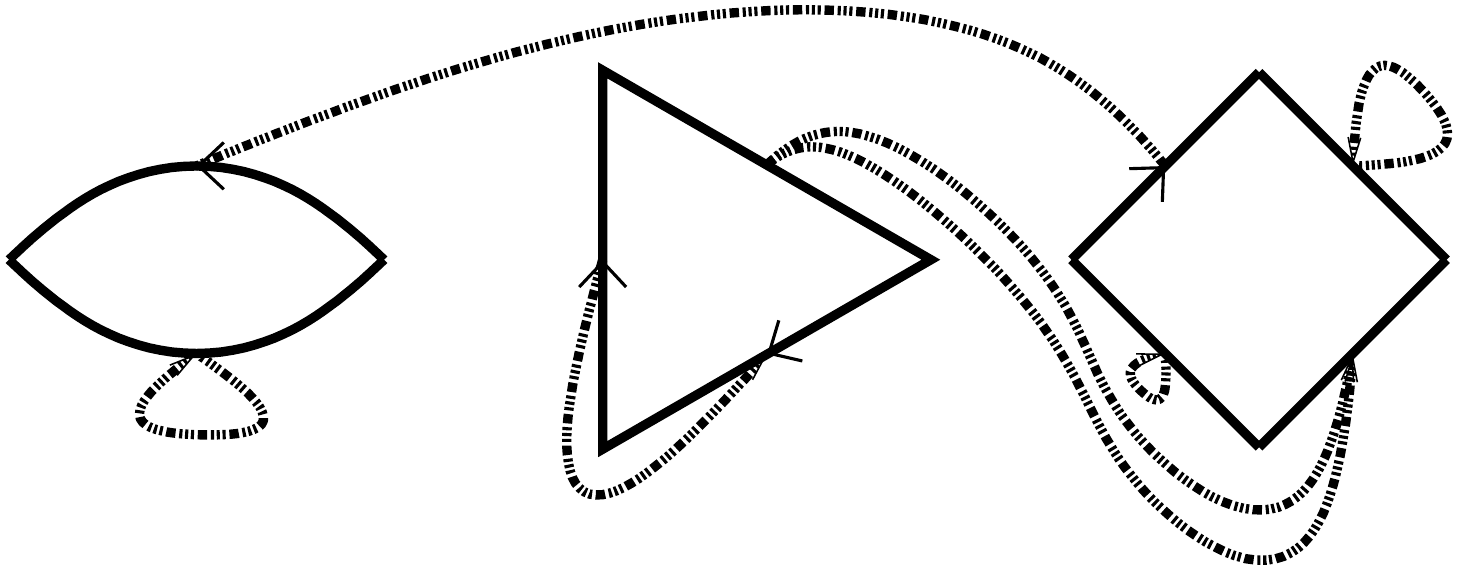}
\caption{A term with some connected terms and some disconnected terms.}
\label{figure: PIE}
\end{figure}

For example, consider the following cumulant of traces of products of centred terms:
\begin{multline*}
k_{3}\left\{\mathrm{tr}\left[\left(Z-\mathbb{E}\left(\mathrm{tr}\left(Z\right)\right)\right)\left(W_{2}^{\left(\lambda_{2}\right)}-\mathbb{E}\left(\mathrm{tr}\left(W_{2}^{\left(\lambda_{2}\right)}\right)\right)\right)\right],\right.\\\left.\mathrm{tr}\left[\left(W_{1}^{\left(\lambda_{3}\right)}-\mathbb{E}\left(\mathrm{tr}\left(W_{1}^{\left(\lambda_{3}\right)}\right)\right)\right)\left(Z^{T}Z^{T}-\mathbb{E}\left(\mathrm{tr}\left(Z^{T}Z^{T}\right)\right)\right)\right],\right.\\\left.\mathrm{tr}\left[\left(W_{2}^{\left(\lambda_{6}\right)}-\mathbb{E}\left(\mathrm{tr}\left(W_{2}^{\left(\lambda_{6}\right)}\right)\right)\right)\left(Z^{T}-\mathbb{E}\left(\mathrm{tr}\left(Z^{T}\right)\right)\right)\right.\right.\\\left.\left.\left(W_{2}^{\left(\lambda_{8}\right)}-\mathbb{E}\left(\mathrm{tr}\left(W_{2}^{\left(\lambda_{8}\right)}\right)\right)\right)\left(W_{1}^{\left(\lambda_{9}\right)}-\mathbb{E}\left(\mathrm{tr}\left(W_{1}^{\left(\lambda_{9}\right)}\right)\right)\right)\right]\right\}\textrm{.}
\end{multline*}
(Note that we have decided to consider the two Ginibre matrices in the second trace as a single term, so we have centred their product.)  In Figure~\ref{figure: PIE}, the second term in the first trace \(W_{2}^{\left(\lambda_{2}\right)}\), the second term in the second trace \(Z^{T}Z^{T}\), and the first and third terms in the third trace \(W_{2}^{\left(6\right)}\) and \(W_{2}^{\left(8\right)}\) are all disconnected.  (Again, even though the two Ginibre matrices in the second trace are connected to each other, this term is not connected to any other term, so it is considered disconnected.)  The set of indices of disconnected terms is then \(\left\{2,4,5,7\right\}\subseteq\left[8\right]\)).  It then appears with sign \(\left(-1\right)^{\left|K\right|}\) in \(g\left(K\right)\) for any \(K\) with \(K\subseteq\left\{2,4,5,7\right\}\), but only appears in \(f\left(\left\{2,4,5,7\right\}\right)\).

We can think of the alternating sign \(\left(-1\right)^{\left|K\right|}\) as subtracting out all terms with one disconnected term (with over-counting), then compensating by adding back those with two disconnected terms (again with over-counting), then again compensating by subtracting those with three disconnected terms, and so on.  We are left with only those terms where no term is disconnected.

\section{Genus Expansions and Asymptotics}
\label{section: genus}

\begin{figure}
\centering
\scalebox{0.5}{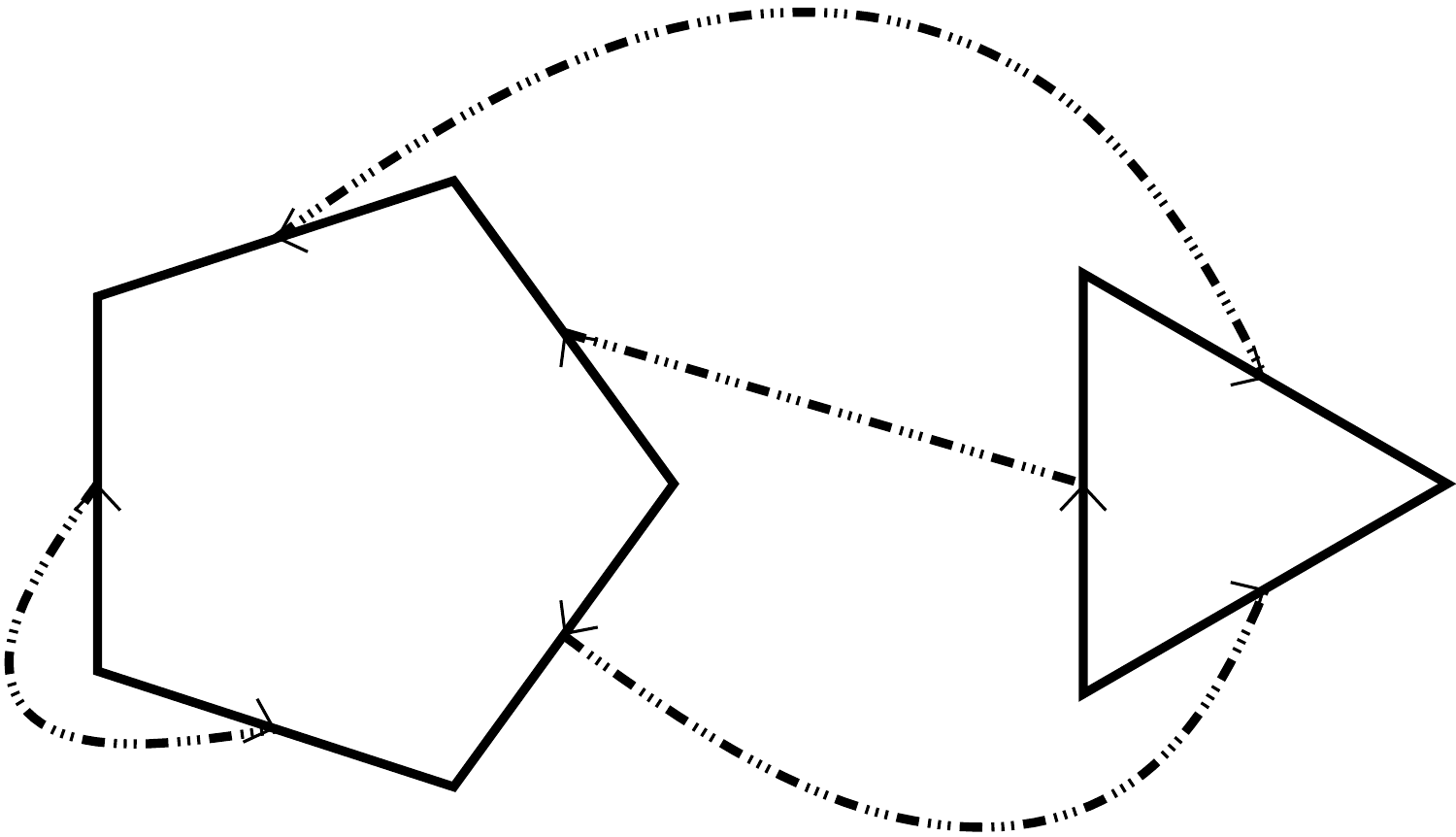}
\caption{An example of a sphere constructed from the faces shown in Figure~\ref{figure: faces}.  The face on the right has been flipped over, so it is the back that is visible.}
\label{figure: noncrossing}
\end{figure}

Moments and cumulants of each of the matrix models we consider can be put in a form similar to what is called a genus expansion, or more accurately in this case an Euler characteristic expansion.  Such a form is a sum of terms corresponding to surfaces, where the order of a term in the size of the matrix \(N\) depends only on the Euler characteristic of the surface.  Highest order terms correspond to highest Euler characteristic, either spheres or collections of spheres.  Gluing diagrams for spheres or collections of spheres can also be drawn as noncrossing diagrams, since crossings require handles to accommodate them.  In the real case, in which there are often two possible directions in which edges may be identified, the edge identifications must not be twisted, since the nonorientable surfaces have lower Euler characterstic than the sphere.  However, we must now consider all relative orientations of the faces or circles.  Figure~\ref{figure: noncrossing} shows an example of a noncrossing diagram on the faces from Figure~\ref{figure: faces} which results in a sphere.  However, the triangle is flipped over (that is, its orientation relative to that of the pentagon is reversed) in order to demonstrate that the diagram is noncrossing.  It is in this orientation that the edge-identifications are untwisted.

The classification theorem for compact connected surfaces states that all connected, compact surfaces are either spheres (Euler characteristic \(2\)), connected sums of \(n\) tori (Euler characteristic \(2-2n\)), or connected sums of \(n\) projective planes (Euler characteristic \(2-n\)) (see, for example, \cite{MR1280460}, section 21, for a combinatorial proof).  The highest order terms of cumulants, which correspond to a single sphere, can be drawn as connected noncrossing diagrams.  In particular, covariances correspond to connected noncrossing diagrams on two circles, or annular noncrossing diagrams, which are studied in \cite{MR0142470, MR2052516}.

\begin{figure}
\centering
\scalebox{0.75}{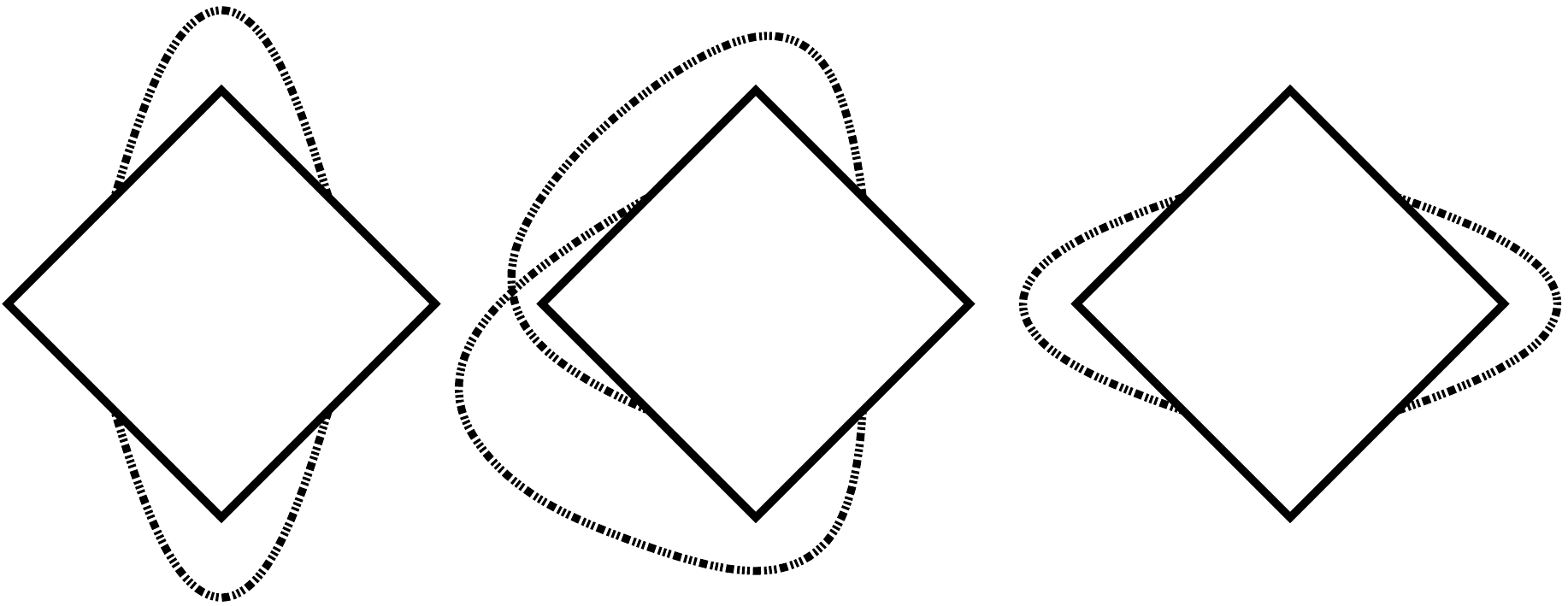}
\caption[Pairings on the edges of a square face]{The three possible pairings on the edges of a square face; the ones on the left and right are noncrossing.}
\label{figure: Ginibre pairings}
\end{figure}

\begin{example}
There are three pairings on \(4\) points, shown in Figure~\ref{figure: Ginibre pairings}.  Two are noncrossing (the one on the left and the one on the right).  To contribute a highest-order term, the edge-identifcations must also be untwisted.  If the matrices are Ginibre matrices, it must therefore be a \(*\)-pairing.  For the \(4\)-gon on the left in Figure~\ref{figure: Ginibre matrices}, only the centre and right pairings are \(*\)-pairings, so only one term contributes to the asymptotic moment.  We get:
\[\lim_{N\rightarrow\infty}\mathbb{E}\left(\mathrm{tr}\left(ZY_{1}ZY_{2}Z^{T}Y_{3}Z^{T}Y_{4}\right)\right)=\mathbb{E}\left(\mathrm{tr}\left(Y_{1}Y_{3}\right)\mathrm{tr}\left(Y_{2}\right)\mathrm{tr}\left(Y_{4}\right)\right)\textrm{,}\]
or, if \(Y_{k}=I_{N}\) for all \(k\), we can calculate by counting diagrams:
\[\lim_{N\rightarrow\infty}\mathbb{E}\left(\mathrm{tr}\left(ZZZ^{T}Z^{T}\right)\right)=1\textrm{.}\]
For the \(4\)-gon on the left in Figure~\ref{figure: Ginibre matrices}, the left and right pairings (but not the centre one) are \(*\)-pairings, so two terms contribute to the asymptotic moment.  We get:
\begin{multline*}
\lim_{N\rightarrow\infty}\mathbb{E}\left(\mathrm{tr}\left(ZY_{5}Z^{T}Y_{6}ZY_{7}Z^{T}Y_{8}\right)\right)\\=\mathbb{E}\left(\mathrm{tr}\left(Y_{1}\right)\mathrm{tr}\left(Y_{2}Y_{4}\right)\mathrm{tr}\left(Y_{3}\right)\right)+\mathbb{E}\left(\mathrm{tr}\left(Y_{1}Y_{3}\right)\mathrm{tr}\left(Y_{2}\right)\mathrm{tr}\left(Y_{4}\right)\right)
\end{multline*}
or
\[\lim_{N\rightarrow\infty}\mathbb{E}\left(\mathrm{tr}\left(ZZ^{T}ZZ^{T}\right)\right)=2\textrm{.}\]
\end{example}

The asymptotic values of the fluctuations are given by the two-face maps appropriate to the matrix model which form a single sphere.  These correspond to the connected noncrossing permutations that can be drawn on the two faces.  In the real case, we must consider whether a permutation is noncrossing on both possible relative orientations of the two faces.

\begin{figure}
\centering
\scalebox{0.75}{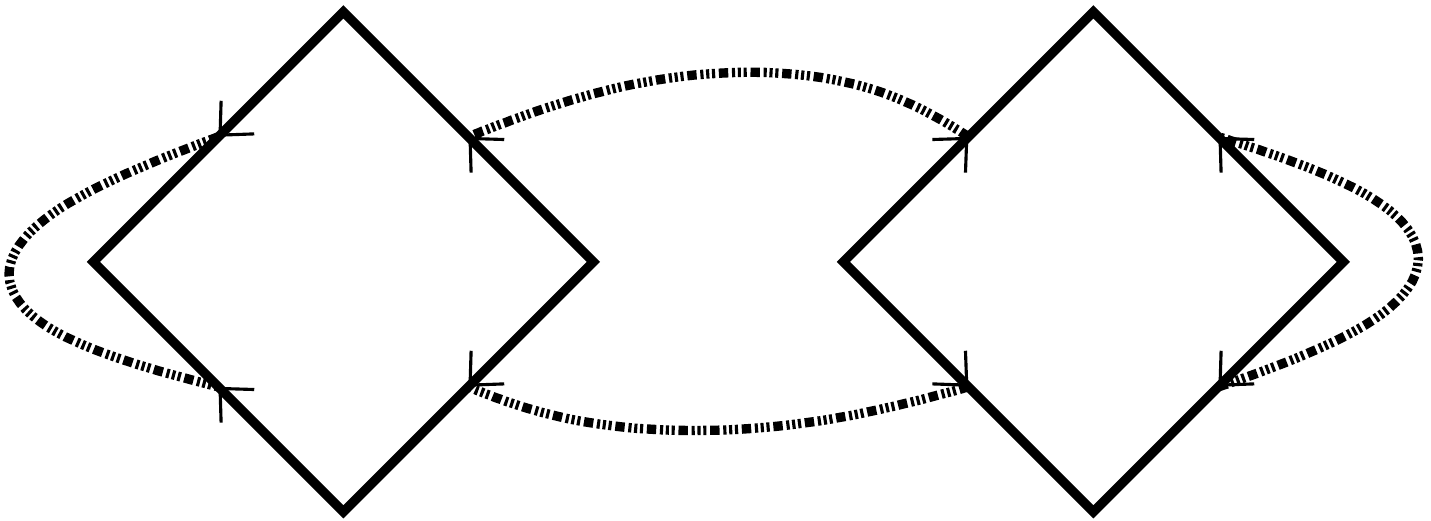}
\caption[Ginibre matrix example with relative orientation of the faces reversed]{The faces and pairing from Figure~\ref{figure: Ginibre matrices} with the face on the right flipped (over a horizontal axis).}
\label{figure: Ginibre fluctuations}
\end{figure}

The fluctuations of the Ginibre matrices are given by the number of annular noncrossing \(*\)-pairings under both relative orientations of the two faces, where the unstarred terms become starred and {\em vice versa} when the orientation of a face is reversed (compare Figure~\ref{figure: Ginibre fluctuations} with Figure~\ref{figure: Ginibre matrices}).

\begin{example}
\label{example: Ginibre fluctuations}
To calculate the asymptotic covariance of \(\mathrm{tr}\left(ZZZ^{T}Z^{T}\right)\) and \(\mathrm{tr}\left(ZZ^{T}ZZ^{T}\right)\), we construct the faces shown in Figure~\ref{figure: Ginibre matrices}.  We show the same faces with the same pairing in Figure~\ref{figure: Ginibre fluctuations}; however, we have flipped the face on the right over its horizontal axis.  On the ``back'' of the face, the transpose of the matrix on the front is shown.  We can see that this pairing is a \(*\)-pairing under this relative orientation, where it is also noncrossing.

We may calculate the asymptotic covariance by counting the annular noncrossing \(*\)-pairings on both relative orientations of the faces, which we may do by inspection in this case.  The pairing must connect the two faces, so there must be at least two pairings between the faces, which must pair an adjacent pair of edges to another adjacent pair of faces to avoid crossings.  We find that there are \(4\) contributing pairings in the relative orientation shown in Figure~\ref{figure: Ginibre matrices} and \(4\) in the relative orientation shown in Figure~\ref{figure: Ginibre fluctuations}.  Thus,
\[\lim_{N\rightarrow\infty}k_{2}\left(\mathrm{Tr}\left(ZZZ^{T}Z^{T}\right),\mathrm{Tr}\left(ZZ^{T}ZZ^{T}\right)\right)=8\textrm{.}\]
\end{example}

\section{Freeness}
\label{section: freeness}

\begin{figure}
\centering
\scalebox{0.75}{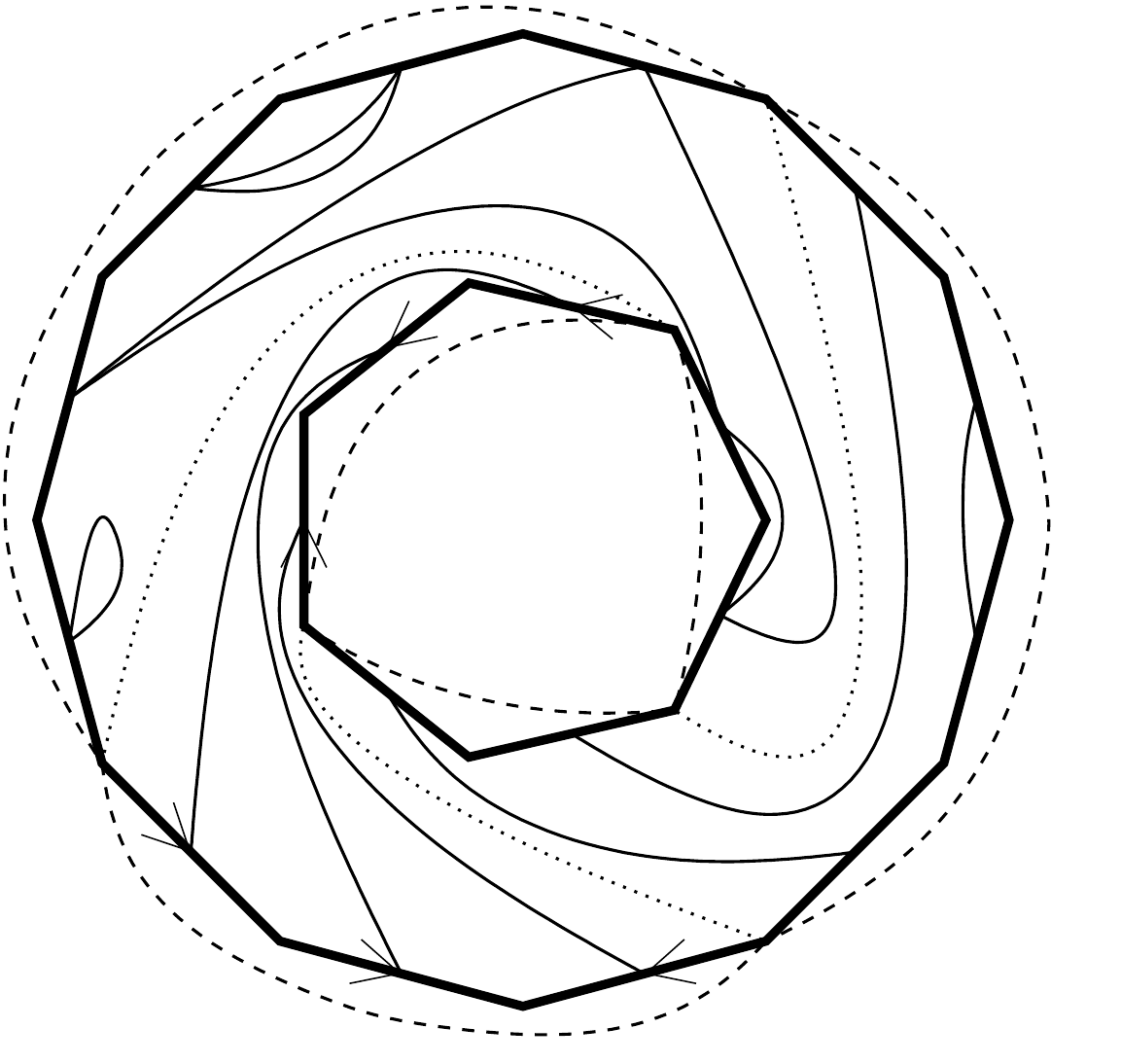}
\caption{An example of a spoke diagram.  Dashed lines mark which edges belong to a single term, and dotted lines separate the spokes.}
\label{figure: spoke diagram}
\end{figure}

In Lemma~6.1 and Theorem~6.2 of the main paper, we consider expressions where the centred terms come from alternating or cyclically alternating algebras, so the terms cannot be connected to their neighbours.  If we wish to consider the terms which survive asymptotically, we require that the diagram be noncrossing.  If a term is connected to a non-neighbouring term on the same face, this separates the terms between the connected terms from all the others.  This eventually forces an unconnected term (forbidden by the centring of the terms) or a term paired with its neighbour (forbidden by the cyclic alternating of the terms).  In the one-trace expressions, we are left with no possible diagrams, so the quantity vanishes asymptotically, meaning that the matrices are asymptotically free.  In the two-trace expressions, a term may instead be connected to a term on the other face.  This significantly restricts the sorts of diagrams which are possible: specifically, the terms must be connected in a sort of spoke arrangement, as in Figure~1 of the main paper (with an example shown in Figure~\ref{figure: spoke diagram} here).  In the real case, both relative orientations of the two faces must be considered.

\begin{figure}
\centering
\scalebox{0.75}{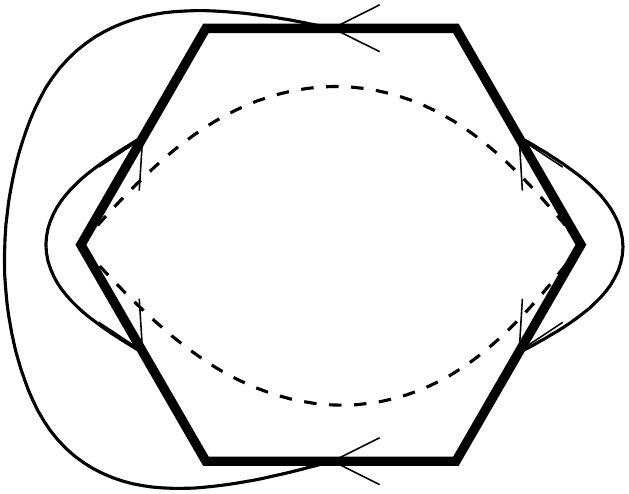}
\caption{The spoke on \(A_{1}\) and \(B_{2}^{T}\), which corresponds to a noncrossing diagram on a disc.}
\label{figure: spoke}
\end{figure}

For a given spoke diagram, the contribution is the sum over all choices of all possible diagrams (those which are connected, noncrossing, untwisted, and compatible with the matrix model) connecting each pair forming a spoke.  Since the contribution of each spoke is multiplicative, we can factor this as the product over the spokes of the sum over the possible diagrams on that spoke.  These diagrams will be disc-noncrossing on the two intervals comprising the spoke (that is, those that would contribute to the limit value of the moment, with an example shown in Figure~\ref{figure: spoke}), but we must subtract out the disconnected diagrams (the product of a choice of a possible term on each interval individually, that is, the product of the moments on each interval individually), so we get the contribution 
\[\mathbb{E}\left(\mathrm{tr}\left(A_{k}B_{l}^{\left(\pm 1\right)}\right)\right)-\mathbb{E}\left(\mathrm{tr}\left(A_{k}\right)\right)\mathbb{E}\left(\mathrm{tr}\left(B_{l}^{\left(\pm 1\right)}\right)\right)\textrm{.}\]

\begin{example}
Let colour \(1\) correspond to a Ginibre matrix \(Z\), colour \(2\) correspond to a GOE matrix \(T\), and colour \(3\) correspond to a Wishart matrix \(W:=X^{T}X\), \(X\) as in the definition.  Let \(A_{1}=ZZZ^{T}\), \(A_{2}=T^{2}\), \(A_{3}=W^{2}\), \(B_{1}=W^{5}\), \(B_{2}=ZZZ^{T}\), and \(B_{3}=T^{4}\).  If we consider the reversed spoke diagram with \(k=1\) (corresponding to the bottom-centre diagram in Figure~1 of the main paper), Figure~\ref{figure: spoke diagram} shows a possible diagram.

By inspection, we can see that this is the only spoke diagram that will accommodate a diagram on \(A_{1}\) and \(B_{2}^{T}\), the two members of the algebra generated by \(Z\).  Again by inspection, we can see that there is only one other diagram on the face shown in Figure~\ref{figure: spoke}, and it also connects \(A_{1}\) to \(B_{2}^{T}\), so the contribution of this spoke is
\[\mathbb{E}\left(\mathrm{tr}\left(ZZZ^{T}\left(ZZZ^{T}\right)^{T}\right)\right)-\mathbb{E}\left(\mathrm{tr}\left(ZZZ^{T}\right)\right)\mathbb{E}\left(\mathrm{tr}\left(\left(ZZZ^{T}\right)^{T}\right)\right)=2\textrm{.}\]
Using the results cited in Remark~5.13 of the main paper, we may calculate that:
\[\mathbb{E}\left(\mathrm{tr}\left(T^{6}\right)\right)-\mathbb{E}\left(\mathrm{tr}\left(T^{2}\right)\right)\mathbb{E}\left(\mathrm{tr}\left(T^{4}\right)\right)=5-1\cdot 2=3\]
and
\[\mathbb{E}\left(\mathrm{tr}\left(W^{7}\right)\right)-\mathbb{E}\left(\mathrm{tr}\left(W^{2}\right)\right)\mathbb{E}\left(\mathrm{tr}\left(W^{5}\right)\right)=429-2\cdot 42=345\textrm{,}\]
so
\[\lim_{N\rightarrow\infty}k_{2}\left(\mathrm{Tr}\left(\mathaccent"7017{A}_{1}\mathaccent"7017{A}_{2}\mathaccent"7017{A}_{3}\right),\mathrm{Tr}\left(\mathaccent"7017{B}_{1}\mathaccent"7017{B}_{2}\mathaccent"7017{B}_{3}\right)\right)=2070\textrm{.}\]
\end{example}

\bibliography{appendix}
\bibliographystyle{plain}

\end{document}